\documentstyle{article}
\begin{document}
\input{amssym}
\def\be{\begin{eqnarray}}
\def\ee{\end{eqnarray}}
\def\RR{{\Bbb R}}
\def\Span{\mathrm{Span}}
\def\di{\displaystyle }
\def\sgn{\mathrm{sgn}}
\title{Lie symmetries of radiation natural convection flow past an inclined surface}
\author{M.~Nadjafikhah \thanks{School of Mathematics, Iran University
of Science and Technology, Narmak, Tehran, 1684613114, I.R.~IRAN.
e-mail: m\_nadjafikhah@iust.ac.ir} , M.~Abdolsamadi,
E.~Oftadeh\thanks{Faculty of Mathematics,
 North Tehran Branch of Islamic Azad University , Pasdaran, Tehran, I.R.IRAN. e-mail:
m.abdolsamadi@gmail.com, ela.oftadeh@gmail.com} }
%
%
\maketitle
\abstract{Applying Lie symmetry method, we find the most general
Lie point symmetry groups of the radiation natural convection flow
equation(RNC). Looking the adjoint representation of the obtained
symmetry group on its Lie algebra, we will find the preliminary
classification of its group-invariant solutions.}

\medskip \noindent {\bf Keywords:} Lie-point symmetries, invariant, optimal
system of Lie sub-algebras.

\medskip \noindent {\bf Mathematics Subject Classification}: 70G65, 58K70, 34C14.
\section{Introduction}
This paper can be viewed as a continuation of the work of paper
\cite{Sivas} where the authors by using Lie group analysis for the
PDE system corresponding to radiation effects on natural
convection heat transfer past an inclined surface (RNC) could
reduce this system to an ODE one and obtain numerical solutions by
applying Runge-Kutta method.

There are various techniques for finding solutions of differential
equations but most of them are useful for a few classes of
equations and applying these techniques for unknown equations is
impossible. Fortunately symmetries of differential equations
remove this problem and give exact solutions.

Symmetry methods for differential equations were discussed  by
S.Lie at first. One of the most important Lie surveys was finding
relationships between continuous transformation group and their
infinitesimal generators that allows us to reduce invariant
conditions of differential equations under group action, that are
complicated because of non linearity, to linear ones.

The main results of applying Lie symmetry method for differential
equations is finding group invariant solutions. A useful way for
reducing equations is finding any subgroup of the symmetry group
and writing invariant solutions with respect to this subalgebra.
This reduced equation is of fewer variables and is easier to
solve. In fact for many important equations arising from geometry
and physics these invariant solutions are the only ones which can
be studied thoroughly.

This type of equations as (RNC) system plays an important role in
engineering and industrial fields. In this paper Lie symmetry
method is applied to find the most general Lie symmetry group and
optimal system of (RNC). Finally, we gain group-invariant
solutions of the reduced system.
\section{Symmetries of RNC System}
Consider the (RNC) system:
\be &&\frac{\partial u}{\partial x}+\frac{\partial v}{\partial
y}=0,\nonumber\\
&&u\frac {\partial u}{\partial x}+v\frac{\partial u}{\partial
y}=\frac{\partial^2 u}{\partial y^2}+Gr\theta \cos
\alpha,\\
&&u\frac{\partial \theta}{\partial x}+v\frac{\partial
\theta}{\partial
y}=\frac{1}{Pr}(1+4R)\frac{\partial^2\theta}{\partial
y^2}.\nonumber\ee
where  $\di Gr=\frac{g\beta(T_W-T_{\infty})\nu}{U_{\infty}^3}$ is
the Grashof number, $\di Pr=\frac{\rho c_p\nu}{k}$ is the Prandtl
number and $\di R=\frac{4\sigma_0T_{\infty}^3}{3k^8}$ is the
radiation parameter \cite{Sivas}. Now the infinitesimal method is
implied as follow:

The infinitesimal generator $\textbf{X}$ on total space of the
form:
\be  \textbf{X} = \sum_{i=1}^p\xi^i(x,u)\partial_{x^i} +
\sum_{j=1}^q\eta^j(x,u)\partial_{u^j},\ee
has the $n^{th}$ order prolongation $\textrm{Pr}^{(n)}X$
(\cite{Olver95}, Th 4.16). Applying the fundamental infinitesimal
symmetry criterion (\cite{Olver95}, Th 6.5) on $X$ as follow:
\be\label{eq:11}
\textrm{Pr}^{(n)}\textbf{X}\big[\Delta_\nu(x,u^{(n)})\big]=0,\quad\nu=1,\cdots,l,\quad\textrm{whenever}\quad
\Delta_\nu(x,u^{(n)})=0. \ee
obtains the infinitesimal generators of the symmetry group.

The vector field associated with RNC system is of the form:
\be X=\xi_1(x,y,u,v)\partial_ x+\xi_2(x,y,u,v)\partial_y \hspace{35mm}\nonumber\\
+\varphi_1
(x,y,u,v)\partial_u+\varphi_2(x,y,u,v)\partial_v+\varphi_3(x,y,u,v)\partial_t.\hspace{9mm}\ee
Since in RNC system, second order derivatives are appeared,
symmetry generators are obtained by applying (\ref{eq:11}) for the
second prolongation of $X$. The vector field $\textbf{X}$
generates a one parameter symmetry group of RNC
$=(\Delta_1,\Delta_2,\Delta_3,\Delta_4)$ if and only if equations
(\ref{eq:11}) hold for $\nu=1,2,3,4$.
So symmetry group of RNC system is spanned by following
infinitesimal generators:
\begin{eqnarray*}
\begin{array}{ll} \di X_1=\partial_x,\!\!\quad\quad&\!\!\quad\quad
\di X_2=\partial_y,\\[3mm]
\di X_3=x\,\partial_x\!+\!u\partial_u+t\partial_t,
\!\!\quad\quad&\!\!\quad\quad\di
X_4=2x\partial_x\!+\!y\partial_y-v\partial_v-2t\partial_t.\\[3mm]
\end{array}
\end{eqnarray*}
The Lie table of Lie algebra $\goth{g}$ for RNC system is given
below:
\begin{eqnarray*}
\begin{array}{c|cccc}
 \!\!\!\quad&\!\!\!\quad X_1 \!\!\!\quad&\!\!\!\quad X_2 \!\!\!\quad&\!\!\!\quad X_3 \!\!\!\quad&\!\!\!\quad X_4 \\
\hline\
{X}_1 \!\!\!\quad&\!\!\!\quad 0 \!\!\!\quad&\!\!\!\quad  0 \!\!\!\quad&\!\!\!\quad X_1 \!\!\!\quad&\!\!\!\quad 2X_1\\
{X}_2 \!\!\!\quad&\!\!\!\quad  0 \!\!\quad&\!\!\quad  0 \!\!\!\quad&\!\!\!\quad  0 \!\!\!\quad&\!\!\!\quad X_2 \\
{X}_3 \!\!\!\quad&\!\!\!\quad  -X_1 \!\!\!\quad&\!\!\!\quad  0 \!\!\!\quad&\!\!\!\quad  0 \!\!\!\quad&\!\!\!\quad  0\\
{X}_4 \!\!\!\quad&\!\!\!\quad  -2X_1 \!\!\!\quad&\!\!\!\quad  -X_2 \!\!\!\quad&\!\!\!\quad  0 \!\!\!\quad&\!\!\!\quad  0\\
\end{array}
\end{eqnarray*}\hfill\
\paragraph{Theorem 1.}
{\it If $g_k(h)$ be the one parameter group generated by $X_k$,
$k=1,\cdots,4$, then}
\be
\begin{array}{rcccl}
g_1&:&(x,y,u,v,t)&\longmapsto&\di (h+x,y,u,v,t),\\[1mm]
g_2&:&(x,y,u,v,t)&\longmapsto&\di (x,h+y,u,v,t),\\[1mm]
g_3&:&(x,y,u,v,t)&\longmapsto&\di (xe^h,y,ue^h,v,te^h),\\[1mm]
g_4&:&(x,y,u,v,t)&\longmapsto&\di (xe^{2h},ye^h,u,ve^{-h},te^{-2h}).\\[1mm]
\end{array}
\ee
\section{Optimal System}
The aim of Lie theory is classifying the invariant solutions and
reducing equations. Some of Lie algebras contain different
subalgebras, so classifying them plays an important role in
transforming equations into easier ones. Therefore we are
classifying these subalgebras up to adjoint representation and
finding an optimal system of subalgebras instead of finding an
optimal system of subgroups. \\The adjoint action is given by the
Lie series
\be \mathrm{Ad}(\exp(\varepsilon X_i)X_j) =
X_j-\varepsilon[X_i,X_j]+\frac{\varepsilon^2}{2}[X_i,[X_i,X_j]]-\cdots,\ee
where $[X_i,X_j]$ is the commutator for the Lie algebra,
$\varepsilon$ is a parameter, and $i,j=1,\cdots,4$
(\cite{Olver93}, ch 3.3).
\paragraph{Theorem 2.}
{\em A one-dimensional optimal system of RNC results as follow:}
\be
\begin{array}{l} 1)\;\;X_1,\\ 2)\;\;X_2,\\ 3)\;\;X_3,\end{array}
\qquad
\begin{array}{l} 4)\;\;X_1+X_2,\\ 5)\;\;X_2-X_1,\\ 6)\;\;X_2+X_3,\end{array} \qquad
\begin{array}{l} 7)\;\;X_3-X_2,\\ 8)\;\;X_3+X_4,\\
9)\;\;X_4-X_3.\end{array} \nonumber \ee
\noindent {\it Proof:} Consider the symmetry algebra $\goth{g}$ of
the RNC system with adjoint representation demonstrated in the
below table:
\begin{eqnarray*}
\begin{array}{c|cccc}
 \!\!\!\quad&\!\!\!\quad {X}_1 \!\!\!\quad&\!\!\!\quad {X}_2 \!\!\!\quad&\!\!\!\quad  {X}_3 \!\!\!\quad&\!\!\!\quad  {X}_4 \\
\hline\
{X}_1 \!\!\!\quad&\!\!\!\quad  {X}_1+\varepsilon{X}_3+2\varepsilon{X}_4 \!\!\!\quad&\!\!\!\quad  {X}_2 \!\!\!\quad&\!\!\!\quad  {X}_3 \!\!\!\quad&\!\!\!\quad  {X}_4\\
{X}_2 \!\!\!\quad&\!\!\!\quad  {X}_1 \!\!\quad&\!\!\quad  {X}_2+\varepsilon{X}_4 \!\!\!\quad&\!\!\!\quad  {X}_3 \!\!\!\quad&\!\!\!\quad  {X}_4 \\
{X}_3 \!\!\!\quad&\!\!\!\quad  {e}^{-\varepsilon}{X}_1 \!\!\!\quad&\!\!\!\quad  {X}_2 \!\!\!\quad&\!\!\!\quad  {X}_3 \!\!\!\quad&\!\!\!\quad  {X}_4\\
{X}_4 \!\!\!\quad&\!\!\!\quad  {e}^{-2\varepsilon}{X}_1 \!\!\!\quad&\!\!\!\quad  {e}^{-\varepsilon}{X}_2 \!\!\!\quad&\!\!\!\quad  {X}_3 \!\!\!\quad&\!\!\!\quad  {X}_4\\
\end{array}
\end{eqnarray*}
 Let
 \be X=a_1X_1+a_2X_2+a_3X_3+a_4X_4.\ee
 In this stage our goal is to simplify as many of the coefficient
 $a_i$ as possible.

 The following results are obtained from adjoint application on
 $X$:
\begin{itemize}
\item[1.]
Let $a_4\neq0$. Scaling $X$ if necessary, suppose that $a_4=1$.
Then
 \be X'=a_1X_1+a_2X_2+a_3X_3+X_4.\ee
 If we act on $X'$ with Ad(exp($a_3X_3$)), the coefficient of $X_1$ can
 be vanished:
 \be X''=a_2X_2+a_3X_3+X_4.\ee
 By applying Ad(exp($a_4X_4$)) on $X''$, we will obtain $X'''$ as
 follow:
\be X'''=a_3X_3+X_4.\ee
\begin{itemize}
\item[1-a.] If $a_3\neq0$, then the coefficient of $X_3$ can be -1
 or 1. Therefore every one-dimensional subalgebra generated by $X$
 in this case is equivalent to $X_3+X_4$, $X_4-X_3$.
\item[1-b.] If $a_3=0$, then every one-dimensional subalgebra generated by $X$
 in this case is equivalent to $X_4$.
\end{itemize}
\item[2.] The remaining one-dimensional subalgebras are spanned by vector
 fields of the form $X$ with $a_4=0$.
\begin{itemize}
\item[2-a.] If $a_3\neq0$, let $a_3=1$. By acting of Ad(exp($a_3X_3$))
 on $X$, we will have
 \be \hat X=a_2X_2+X_3.\ee
\begin{itemize}
\item[2-a-1.] If $a_2\neq0$, then the coefficient of $X_2$ can
 be -1 or 1. So  every one-dimensional subalgebra generated by $X$
 in this case is equivalent to $X_2+X_3$, $X_3-X_2$.
\item[2-a-2.] If $a_2=0$ then every one-dimensional subalgebra
generated by $X_3$.
\end{itemize}
\item[2-b.] Let $a_4=0,a_3=0$.
\begin{itemize}
\item[2-b-1.] If $a_2\neq0$, then we can make the coefficient of
$X_1$ either -1 or 1 or 0. Every one-dimensional subalgebra
generated by $X$ is equivalent to $X_1+X_2$, $X_2-X_1$, $X_2$.
\item[2-b-2.] If $a_2=0$, then every one-dimensional subalgebra
generated by $X$ is equivalent to $X_1$. \hfill\ $\Box$
\end{itemize}
\end{itemize}
\end{itemize}
\section{Invariant Solutions and Reduction }
Now for finding invariant solutions using characteristic method
\cite{Hydon2000}:
\be Q_\alpha\mid_{u=u(x,t)}\equiv
X[u^\alpha-u^\alpha(x,t)]\mid_{u=u(x,t)} \quad \alpha=1,\cdots,M
\ee
where $M$ is the number of dependant variables.
Now invariant solutions for RNC result as follow:\\
\\
 Consider $X_1=\partial_x$. Applying characteristic method on
$X_1$ obtains: $Q_u=-{\frac {\partial }{\partial x}}u \left( x,y
\right), Q_v=-{\frac {\partial }{\partial x}}v \left( x,y \right),
Q_t=-{\frac {\partial }{\partial x}}t \left( x,y \right)$.
Solutions of this system are of the form: $u={F_1} \left( y
\right), v={F_1} \left( y \right), t={F_1} \left( y \right) $.\\
\\
Consider $X_2=\partial_y$. Applying characteristic method on $X_2$
obtains: $Q_u=-{\frac {\partial }{\partial y}}u \left( x,y
\right), Q_v=-{\frac {\partial }{\partial y}}v \left( x,y \right),
Q_t=-{\frac {\partial }{\partial y}}t \left( x,y \right)$.
Solutions of this system are of the form: $u={F_1} \left( x
\right), v={F_1} \left( x \right), t={F_1} \left( x \right) $.\\
\\
Consider $X_3=x\partial_x+u\partial_u+t\partial_t$. Applying
characteristic method on $X_3$ obtains: $Q_u=u \left( x,y \right)
-x{\frac {\partial }{\partial x}}u \left( x,y
 \right)
, Q_v=-x{\frac {\partial }{\partial x}}v \left( x,y \right), Q_t=t
\left( x,y \right) -x{\frac {\partial }{\partial x}}t \left( x,y
 \right)$.\\
Solutions of this system are of the form: $u=x{F_1}\left( y
\right), v={F_1} \left( y \right), t=x{F_1} \left( y \right) $.\\
\\
Consider $X_4=2x\partial_x+y\partial_y-v\partial_v-2t\partial_t$.
Applying characteristic method on $X_4$ obtains: $Q_u=-2\,x{\frac
{\partial }{\partial x}}u \left( x,y \right) -y{\frac {
\partial }{\partial y}}u \left( x,y \right),
Q_v=-v \left( x,y \right) -2\,x{\frac {\partial }{\partial x}}v
\left( x,y
 \right) -y{\frac {\partial }{\partial y}}v \left( x,y \right)
, Q_t=-2\,t \left( x,y \right) -2\,x{\frac {\partial }{\partial
x}}t \left( x,y \right) -y{\frac {\partial }{\partial y}}t \left(
x,y \right)
$.\\
Solutions of this system are of the form: $u={F_1} \left( {\frac
{y}{\sqrt {x}}} \right), v={F_1} \left( {\frac {y}{\sqrt {x}}}
\right) {\frac {1}{\sqrt {x }}}
, t={F_1} \left( {\frac {y}{\sqrt {x}}} \right) {x}^{-1} $.\\
\\
Consider $X_1+X_2=\partial_x+\partial_y$. Applying characteristic
method on $X_1+X_2$ obtains: $Q_u={\frac {\partial }{\partial x}}u
\left( x,y \right) +{\frac {\partial }{\partial y}}u \left( x,y
\right) , Q_v={\frac {\partial }{\partial x}}v \left( x,y \right)
+{\frac {\partial }{\partial y}}v \left( x,y \right) , Q_t={\frac
{\partial }{\partial x}}t \left( x,y \right) +{\frac {\partial
}{\partial y}}t \left( x,y \right)$. Solutions of this system are
of the form: $u={F_1}\left( -x+y
\right), v={F_1} \left( -x+y \right), t={F_1} \left( -x+y \right) $.\\
\\
Consider $X_2-X_1=\partial_y-\partial_x$. Applying characteristic
method on $X_2-X_1$ obtains: $Q_u={\frac {\partial }{\partial x}}u
\left( x,y \right) -{\frac {\partial }{\partial y}}u \left( x,y
\right) , Q_v={\frac {\partial }{\partial x}}v \left( x,y \right)
-{\frac {\partial }{\partial y}}v \left( x,y \right)
 , Q_t={\frac {\partial }{\partial x}}t \left( x,y \right) -{\frac {\partial
}{\partial y}}t \left( x,y \right)$. Solutions of this system are
of the form: $u={F_1}\left( x+y
\right), v={F_1} \left( x+y \right), t={F_1} \left( x+y \right) $.\\
\\
Consider $X_2+X_3=x\partial_x+\partial_y+u\partial_u+t\partial_t$.
Applying characteristic method on $X_2+X_3$ obtains: $Q_u=u \left(
x,y \right) -x{\frac {\partial }{\partial x}}u \left( x,y
 \right) -{\frac {\partial }{\partial y}}u \left( x,y \right)
 , Q_v=-x{\frac {\partial }{\partial x}}v \left( x,y \right) -{\frac {
\partial }{\partial y}}v \left( x,y \right)
 , Q_t=t \left( x,y \right) -x{\frac {\partial }{\partial x}}t \left( x,y
 \right) -{\frac {\partial }{\partial y}}t \left( x,y \right)
$. Solutions of this system are of the form: $u=x{F_1}\left( -\ln  \left( x \right) +y \right), v={F_1} \left(  -\ln  \left( x \right) +y \right), t=x{F_1} \left( -\ln  \left( x \right) +y \right) $.\\
\\
Consider $X_3-X_2=x\partial_x-\partial_y+u\partial_u+t\partial_t$.
Applying characteristic method on $X_3-X_2$ obtains: $Q_u=u \left(
x,y \right) -x{\frac {\partial }{\partial x}}u \left( x,y
 \right) +{\frac {\partial }{\partial y}}u \left( x,y \right)
 , Q_v=-x{\frac {\partial }{\partial x}}v \left( x,y \right) +{\frac {
\partial }{\partial y}}v \left( x,y \right)
 , Q_t= t \left( x,y \right) -x{\frac {\partial }{\partial x}}t \left( x,y
 \right) +{\frac {\partial }{\partial y}}t \left( x,y \right)
$. Solutions of this system are of the form: $u=x{F_1}\left(  \ln  \left( x \right) +y \right), v={F_1} \left(   \ln  \left( x \right) +y \right), t=x{F_1} \left(  \ln  \left( x \right) +y \right) $.\\
\\
Consider
$X_3+X_4=3x\partial_x+y\partial_y+u\partial_u-v\partial_v-t\partial_t$.
Applying characteristic method on $X_3+X_4$ obtains: $Q_u=u \left(
x,y \right) -3\,x{\frac {\partial }{\partial x}}u \left( x,y
 \right) -y{\frac {\partial }{\partial y}}u \left( x,y \right)
 , Q_v=-v \left( x,y \right) -3\,x{\frac {\partial }{\partial x}}v \left( x,y
 \right) -y{\frac {\partial }{\partial y}}v \left( x,y \right)
 , Q_t=-t \left( x,y \right) -3\,x{\frac {\partial }{\partial x}}t \left( x,y
 \right) -y{\frac {\partial }{\partial y}}t \left( x,y \right)
$. Solutions of this system are of the form: $u={F_1} \left(
{\frac {y}{\sqrt [3]{x}}} \right) \sqrt [3]{x}, v={F_1} \left(
{\frac {y}{\sqrt [3]{x}}} \right) {\frac {1}{\sqrt [3]{x}}}
 ,
t={F_1} \left( {\frac {y}{\sqrt [3]{x}}} \right) {\frac {1}{\sqrt
[3]{x}}}
$.\\
\\
Consider
$X_4-X_3=x\partial_x+y\partial_y-u\partial_u-v\partial_v-3t\partial_t$.
Applying characteristic method on $X_4-X_3$ obtains: $Q_u=-u
\left( x,y \right) -x{\frac {\partial }{\partial x}}u \left( x,y
 \right) -y{\frac {\partial }{\partial y}}u \left( x,y \right)
 , Q_v=-v \left( x,y \right) -x{\frac {\partial }{\partial x}}v \left( x,y
 \right) -y{\frac {\partial }{\partial y}}v \left( x,y \right)
 , Q_t=-3\,t \left( x,y \right) -x{\frac {\partial }{\partial x}}t \left( x,y
 \right) -y{\frac {\partial }{\partial y}}t \left( x,y \right)
$. Solutions of this system are of the form: $u={F_1}  \left(
{\frac {y}{x}} \right) {x}^{-1}, v={F_1} \left( {\frac {y}{x}}
\right) {x}^{-1}, t={F_1} \left( {\frac {y}{x}} \right) {x}^{-3}
$.\\

\section{Conclusion}
In this paper by applying infinitesimal symmetry methods, we find
optimal system and finally could reduce the RNC system and find
invariant solutions.


\begin{thebibliography}{}
%
\bibitem{Bluman}  {\sc G.W. Bluman} and {\sc S. Kumei}, {\em Symmetries and Differential Equations},
Springer, New York, 1989.
%
\bibitem{Hydon2000}  {\sc P. Hydon} , {\em Symmetries Methods for Differential Equations, A Beginner's Guid}, Cambridge: Cambridge University Press.
%
\bibitem{Nadj2010}  {\sc M. Nadjafikhah} , {\em Classification of similarity solutions for inviscid burgers'
equation}, Adv. appl. Clifford alg. 20 (2010), 71-77.
%
\bibitem{Nadj2011} {\sc M. Nadjafikhah}, {\em Lie Symmetries of Inviscid Burgers'
Equation}, Adv. Appl. Clifford Alg., 19. (2009) 101-112.
%
\bibitem{Olver93} {\sc P.J. Olver}, {\em Applications of Lie
Groups to Differential Equations}, Springer, New York, 1993.
%
\bibitem{Olver95} {\sc P.J. Olver}, {\em Equivalence, Invariants,
and Symmetry}, Cambridge University Press, 1995.
%
\bibitem{Sivas} {\sc S. Sivasankaran}, {\sc M. Bhuvaneswari}, {\sc
P. Kandaswamy} and {\sc E.K. Ramasami}, {\em Lie group analysis of
radiation natural convection flow past an inclined surface},
Communications in Nonlinear Science and Numerical Simulation 13
(2008) 269-276.
%
\end{thebibliography}
\end{document}